\documentclass[a4paper]{amsart}
\usepackage{amssymb,latexsym}

\begin{document}

\title[Commutators of Commutators as Cubes]%
{Writing Commutators of Commutators\\as Products of Cubes}

\author[C. Ramsay]{Colin Ramsay}

\address{School of Information Technology and Electrical Engineering\\
The University of Queensland, Queensland 4072, Australia}

\email{uqcramsa@uq.edu.au}

\date{July 29, 2015}



\begin{abstract}
It is known that commutators of commutators can be written as products of cubes,
  with the current upper bound on the number of cubes being $60$.
We discuss how proofs extracted via coset enumeration can be used to
  investigate this problem, and exhibit a rewriting using only $14$ cubes.
\end{abstract}

\maketitle

\section{Introduction}

For elements $x_1$, $x_2$, $x_3$, $x_4$ of a group, the commutator 
  $[x_1,x_2]$ is $x_1^{-1}x_2^{-1}x_1x_2$ and the
  commutator of commutators is $[[x_1,x_2],[x_3,x_4]]$. 
In \cite{Lyndon1980} Lyndon notes that every commutator of commutators is a 
  product of a bounded number of cubes, and asks for an explicit formula with 
  the number of cubes as small as possible.
Hegarty \cite{Hegarty1996} established an upper bound of 85 cubes,
  while Akhavan-Malayeri \cite{Ak-Mal2009} improved the bound to 60.

The Todd-Coxeter coset enumeration process \cite{Sim,TC} is a
  method for enumerating the cosets of a subgroup of a finitely presented group.
The workings involved in an enumeration are usually discarded, but
  if they are recorded they can be used to extract proofs of relations in the
  group \cite{HR5,Lee2}.
These proofs are in the form of \emph{proof-words}, which express a word of 
  interest as a product of subgroup generators and conjugates of relators.
This technique was used by Havas \cite{Hav3} to express the fifth Engel word 
  as a product of $250$ fourth powers, and we use a similar method to express
  the commutator of commutators as a product of $14$ cubes.

\section{Background}

Let $F_4$ be the free group of rank $4$ generated by $\{x,y,z,w\}$, and set
  $C = [[x,y],[z,w]]$.
We use upper-case to represent generator inverses 
  and typewriter font when writing proof-words.
So, for example, $\verb|x|=x$ and $\verb|X| = X = x^{-1}$\kern-2pt.
When manipulating proof-words the following properties of conjugation
   and powers are useful.
Let $u,v \in F_4$ and $n>0$.
Then $v^{-1}u^nv = (u^n)^v = (u^v)^n$ and $u^{-n} = (u^{-1})^n$\kern-2pt.
  So we can write $u^{nv}$ and $u^{-nv}$ without any ambiguity, 
  and $u^nv = vu^{nv}$\kern-2pt.

Groups of exponent three are metabelian \cite[Lemma~5.11]{MKS66}, and thus
  commutators of commutators are trivial is such groups.
Recall that the free Burnside group $B(r,3)$ is the largest group of exponent
  three with $r$ generators.
These Burnside groups are finite and have orders $3^s$\kern-2pt,
  where $s = r(r^2+5)/6$ \cite[Theorem~5.24]{MKS66}.

Fix $G = B(4,3)$.
Given a presentation for $G$, $C$ is in the normal closure of the relators, and 
  proofs of $C$'s triviality extracted via coset enumeration consist of products
  of cubes.
For non-trivial subgroups proof-words will generally include subgroup 
  generators.
However these define a trivial word in $G$ and can be replaced by a product of
  cubes.

The utility {\sc peace} \cite{HR5} was used to perform coset enumerations and
  proof extractions.
The proofs depend on the presentation and subgroup used, the enumerator
  settings, and the internal details of the proof extraction 
  process \cite{Hav3,HR5}.
To obtain good results -- in our case, a proof with a small number of
  cubes -- a large number of runs with varying parameters is generally
  required.

We start with a list of all \emph{base-words} in the alphabet
  $\{x,X,y,Y,z,Z,w,W\}$ for lengths, say, 1--8.
A base-word is a freely and cyclically reduced word which is distinct
  from all the other words under inversion and cyclic permutation.
Relators are third powers of base-words, and the cubes in our 
  expressions for $C$ will be conjugates of these relators.

The exact form of a presentation is a significant source of variation in
  proof-words, so we build our presentations using random selections of 
  base-words, with each word randomly inverted or cyclically permuted.
We limited our attention to presentations which defined $G$
  and enumerated efficiently (that is, the total number of cosets defined did
  not greatly exceed the index).

The choice of subgroup generators also has a significant effect on the 
  proof-words generated.
We experimented with a variety of subgroup generators, but most of our runs were
  done using some proper subset of the generators of $F_4$.
That is, the subgroups were $B(r,3)$, for $1 \leqslant r \leqslant 3$.

Another useful source of variation in proof-words is the choice of 
  word to be proved.
Instead of extracting a proof for $C$ in each run
  we extract $16$ proofs, one for each of $C$\kern-1pt, 
  $C^{W}$\kern-3pt, $C^{W\!Z}$\kern-3pt, etc.
These proof-words need not be conjugates of each other, and are typically
  variable in structure and length.
Any proofs of interest can conjugated to prove $C$ whilst leaving the
  number of cubes unchanged.

The process described is in no sense exhaustive, and our results give no
  guidance as to Lyndon's ``as small as possible''\kern-2pt,
  except to reduce the upper bound to $14$.
Note that at least two cubes are required since, in a free group, no 
  non-trivial commutator is a proper power \cite[page~94]{Lyndon1980}.

\section{Results}

The following proof-word was extracted from an enumeration with base-words of
  length one excluded and using the subgroup $\langle x,y,z \rangle$ of $G$.
\begin{verbatim}
 Wzy(YwYwYw)Wy(WXyWXyWXy)YwYZw[Y][Z](zyWzyWzyW)[X][z][y]
 (YZxwYZxwYZxw)[x][Z][x](XzXWXzXWXzXW)[z][X](xZwxZwxZw)WzX
 (WXWXWX)xw(xwxxwxxwx)WZw[z][X][Y][x][y][Z]
\end{verbatim}
The substrings within parentheses are relators, which may be inverted or cycled,
  while the substrings within brackets are subgroup generators or their
  inverses.
The remaining group generators conjugate the relators.
Note that matched generator/inverse pairs may conjugate multiple relators, and
  that there is no conjugation of subgroup generators.
The entire string freely reduces to $C^{W\!Zw}$\kern-2pt.
After deleting the relators, which are trivial in $G$, all the conjugation
  cancels, leaving only the subgroup element
  $\beta = Y\!ZX\!zyxZxzX\!zXY\!xyZ$.

So this proof-word proves that 
  $C^{W\!Zw} = \beta \in \langle x,y,z \rangle$.
Now $\beta$ is trivial in $G$ and can be rewritten as a product of cubes.
First collect the subgroup generators at the end of the proof-word  by
  introducing appropriate trivial subwords in the tail of the proof-word, just
  before any subgroup generators.
The preceding subgroup generators can now be shuffled to the end of the 
  proof-word by moving the brackets and changing the original subgroup
  generators to conjugation, as illustrated below.
\begin{verbatim}
 ...XzXW)[z][X](xZwxZwxZw)...(xwxxwxxwx)WZw        [z][X][Y]...
 ...XzXW)[z][X](xZwxZwxZw)...(xwxxwxxwx)WZw  xZzX  [z][X][Y]...
 ...XzXW) z  X (xZwxZwxZw)...(xwxxwxxwx)WZwxZ[z][X][z][X][Y]...
\end{verbatim}

To rewrite $\beta$ as a product of cubes we simply need to extract a proof of 
  its triviality using {\sc peace} with a presentation of 
  $\langle x,y,z \rangle = B(3,3)$.
The shortest proof-word for $\beta$ found contained six cubes.
Substituting this for $\beta$ and absorbing some conjugation into the 
  relators yields the following proof-word.
\begin{verbatim}
 Wzy(YwYwYw)W(yWXyWXyWX)wYZwYZ(zyWzyWzyW)Xzy(YZxwYZxwYZxw)xZx 
 (XzXWXzXWXzXW)zX(xZwxZwxZw)Wz(XWXWXW)w(xwxxwxxwx)WZw
 xZXzXYZxzy(YZZYZZYZZ)zY(yzyzyz)ZYxZY(yzXzyzXzyzXz)(ZxZxZx)
 (zXzXzX)yzXyz(ZYxZYxZYx)yZ
\end{verbatim}

If we set $\gamma = W\!zwY\!XyxW\!ZwxZX\!zXY\!Zxzy$, then 
  $C^{W\!Zw} = \gamma\beta$, and the proof-word can be split into two smaller
  proof-words, using the first eight and last six cubes.
A second proof of $C^{W\!Zw}$ with a similar structure was also found.
Set
\begin{align*}
  \delta   &= W\!zwY\!XyxW\!ZwX\!ZY\!zxyZyzY\!,\\
  \epsilon &= yZY\!zY\!X\!ZyzxzXY\!xyZ.    
\end{align*}
Then $C^{W\!Zw} = \delta\epsilon$ and a proof-word for this is given
  below.
Note that the freely-trivial word $zY\!yZ$ has to be introduced immediately
  before the ninth relator to allow the split into valid proof-words.
\begin{verbatim}
 WzwY(XywXywXyw)yW(wywywy)YW(YWYYWYYWY)wyZw(WzYWzYWzY)yZYz
 (ZywyZywyZywy)YXZy(YzxWYzxWYzxW)z(ZwXZwXZwX)xWz(xWxWxW)Zw 
 XZYzxyZy(YzYzYz)(ZyZXZyZXZyZX)xzY(zxzzxzzxz)(ZXZXZX)y
 (YxzYxzYxz)ZXz(ZyZyZy)Y
\end{verbatim}

These proof-words express $C^{W\!Zw}$ as a product of $14$ relators 
  (with total length 120) with 24 \& 21 conjugating generator/inverse pairs
  respectively.
Explicit products of cubes are easily recovered by distributing the conjugation,
  with the first two cubes of the second proof-word being
  $(X\!yw)^{3yW\!Zw}$ and $(wy)^{3Zw}$\kern-2pt.

While we did not find any proof-words with fewer relators, we did find
  ones with shorter overall length.
The following proof-word expresses $C$ as a product of $15$ relators (with total
  length $96$) and $15$ conjugating generator/inverse pairs.
\begin{verbatim}
 y(YYY)(yXyXyX)xY(xxx)XWZw(zzz)Z(ZXZXZX)(xzWxzWxzW)zWz(ZwZwZw)
 w(WWW)x(XwXwXw)Wx(ZXwyZXwyZXwy)XwyX(xYWxYWxYW)(wywywy)YWZw
 (WzYWzYWzY)yZ(xzYxzYxzY)zY(yZyZyZ)Wzw
\end{verbatim}

Note that each of the base-words corresponds to one of the $15$ non-empty
  subsets of the set of generators $\{x,y,z,w\}$, up to generator inversions
  and ordering.
For example, $\{y\}$ and $\{y,x\}$ correspond to $Y$ and $yX$ respectively.

\end{document}